\documentclass[12pt]{amsart}
\usepackage{amsfonts,amsmath,amsthm,amscd,amssymb,latexsym,cite,verbatim,texdraw,floatflt,
caption2}
\RequirePackage{hyperref}

\usepackage[a4paper]{geometry}

\theoremstyle
{plain}

\begin{document}

\title{On balanced coronas of  groups }

\author{Igor Protasov}

\maketitle
\vskip 5pt

{\bf Abstract.}
Let $G$ be an  infinite group, $\kappa$  be an  infinite cardinal,
$\kappa\leq \mid G\mid$  and let  $\mathcal{E}_{\kappa}$    denotes
 a  coarse structure on  $G$ with  the base
 $\{\{ (x,y): y\in F x F\}: F\in [G]^{<\kappa}\}$. We prove that if
    either  $\kappa< \mid G\mid$
     or $\kappa= \mid G\mid$
      and  $\kappa$  is singular then the Higson's corona
      $\nu _{\kappa} (G)$  of the coarse space $(G, \mathcal{E}_{\kappa})$ is a singleton.
      If  $\kappa= \mid G\mid$   and $\kappa$ is regular then $\nu _{\kappa} (G)$
       contains a copy  of the space $U_{\kappa}$ of  $\kappa$-uniform  ultrafilters on $\kappa$.

\vskip 10pt

{\bf MSC: } 54E35,  20F69.
\vskip 10pt

{\bf Keywords:}  coarse structure, slowly oscillating functions, balanced corona.

\section{Introduction}

Let $X$  be a set. A family $\mathcal{E}$ of subsets of $X\times X$ is called a {\it coarse structure } if
\vskip 7pt

\begin{itemize}
\item{}   each $E\in \mathcal{E}$  contains the diagonal  $\bigtriangleup _{X}$,
$\bigtriangleup _{X}= \{(x,x): x\in X\}$;
\vskip 5pt

\item{}  if  $E$, $E^{\prime} \in \mathcal{E}$ then $E\circ E^{\prime}\in\mathcal{E}$ and
$E^{-1}\in \mathcal{E}$,   where    $E\circ E^{\prime}=\{(x,y): \exists z((x,z) \in  E,  \   \ (z, y)\in E^{\prime})\}$,   $E^{-1}=\{(y,x): (x,y)\in E\}$;
\vskip 5pt

\item{} if $E\in\mathcal{E}$ and $\bigtriangleup_{X}\subseteq E^{\prime}\subseteq E  $   then
$E^{\prime}\in \mathcal{E}$;
\vskip 5pt

\item{}  for any   $x,y\in X$, there exists $E\in \mathcal{E}$   such that $(x,y)\in E$.

\end{itemize}
\vskip 7pt

A subset $\mathcal{E}^{\prime} \subseteq \mathcal{E}$  is called a
{\it base} for $\mathcal{E}$  if, for every $E\in \mathcal{E}$, there exists
  $E^{\prime}\in \mathcal{E}^{\prime}$  such  that
  $E\subseteq E ^{\prime}$.
For $x\in X$,  $A\subseteq  X$  and
$E\in \mathcal{E}$, we denote
$E[x]= \{y\in X: (x,y) \in E\}$,
 $E [A] = \cup_{a\in A}   \   \   E[a]$
 and say that  $E[x]$
  and $E[A]$
   are {\it balls of radius $E$
   around} $x$  and $A$.
 We say that  a subset  $B$  of $X$  is {\it bounded}  if there exist
    $E \in  \mathcal{E}$  and $a\in X$  such that  $B\subseteq  E[a]$.

The pair $(X,\mathcal{E})$ is called a {\it coarse space} \cite{b14}  or a {\it ballean } \cite{b10}, \cite{b13}.

Let  $(X, \mathcal{E})$, $(X^{\prime}, \mathcal{E}^{\prime})$  be balleans.
A mapping
$f: X \longrightarrow X^{\prime}$  is called {\it macro uniform} if, for every $E \in \mathcal{E}$, there exists
 $E^{\prime} \in \mathcal{E}^{\prime}$   such that
   $f(E[x])\subseteq  E^{\prime}[f(x)]$ for each $x\in X$.
If $f$   is a bijection such that  $f$  and   $f ^{-1}$  are macro-uniform  then $f$  is called an {\it asymorphism}.

Given a ballean $(X, \mathcal{E})$,  we endow  $X$ with  the discrete topology and take  the points of
the Stone-$\check{C}$ech compactification  $\beta X$ of   $X$   to be the ultrafilters  on $X$   with the points of $X$     identified  with the principal ultrafilters.
For every  subset $A$ of $X$,  we put
$\bar{A}= \{q \in \beta X:  A\in q\}$.
The topology  of  $\beta X$  can be defined  by stating that the family
 $\{\bar{A}:   A\subseteq X\}$
  is a base  for the open sets.

If  $Y$  is a compact Hausdorff space and $f: X\longrightarrow   Y$  then  $f^{\beta}$
 denotes the extension of  $f$  to $\beta X$.

We denote $X^{\sharp} = \{ q \in \beta X: $   each  member  of  $q$  is unbounded in
 $(X, \mathcal{E})\}$ and,  given any  $r,q\in X^{\sharp}$, write  $r\| q$ if there  exists
  $E\in \mathcal{E}$  such  that  $E[R]\in  q$
   for each $R\in r$.
By [6, Lemma 4.1],
$\|$
 is an equivalence on $X^{\sharp}$.
We denote by $\sim$  the  minimal (by inclusion) closed
(in $X^{\sharp}\times X^{\sharp}$) equivalence on  $X$  such  that
$\| \subseteq \sim$.
By [8, Proposition 1],  $r\sim  q$  if and only if
$h^{\beta}(p)= h^{\beta}(q)$ for every  slowly  oscillating  function
 $h: X\longrightarrow [0,1]$.
We recall that a function $f :  (X, \mathcal{E})\longrightarrow  \mathbb{R} $
is {\it slowly oscillating}  if, for every  $\varepsilon > 0$,  there exists a  bounded  subset $B$  of  $X$  such that
$diam \  f(E[x])<  \varepsilon$  for each  $x\in X\setminus B$.

The quotient  $X^{\sharp} / \sim$  is called the  {\it Higson's  corona} \cite{b14}  or {\it corona}  \cite{b6}, \cite{b8}  of the ballean  $(X, \mathcal{E})$.

Now let $G$  be an infinite group, $\kappa$ be an infinite cardinal,  $\kappa\leq \mid G \mid$ and let $\mathcal{L}_{\kappa}$, $\mathcal{R}_{\kappa}$ and $\mathcal{E}_{\kappa}$ be coarse structures on $G$  with the bases
$$\{\{ (x,y): y\in x F\}:  F\in [G]^{<\kappa}\}, $$
 $$\{\{ (x,y): y\in F x\}:  F\in [G]^{<\kappa}\}, $$
  $$\{\{ (x,y): y\in F x F\}:  F\in [G]^{<\kappa}\} $$
and say that $(G, \mathcal{L}_{\kappa})$, $(G, \mathcal{R}_{\kappa})$  and
$(G, \mathcal{E}_{\kappa})$ are {\it left $\kappa$-ballean,  right $\kappa$-ballean} and
{\it balanced  $\kappa$-ballean}  of $G$.
We observe that the mapping  $G\longrightarrow G$, $x\longmapsto x^{-1}$
is an asymorphism between  $(G, \mathcal{L}_{\kappa})$ and $(G, \mathcal{R}_{\kappa})$.

We note that left or right  $\omega$-balleans  play important role in  {\it  Geometric Group Theory, }  see [5, Chapter 4]  and  [3, Chapter 3].
Balanced   $\omega$-balleans  were introduce in [10, Chapter 12].
For  right and balanced  $\kappa$-balleans   see \cite{b12},  \cite{b1}  and \cite{b11}  respectively.

\section{Results}

We  say that a function  $h: (G,  \mathcal{E}_{\kappa} )\longrightarrow  \mathbb{R}$
 is  {\it convergent  at infinity}  if there  exists  $r\in \mathbb{R}$  such that,
  for every  $\varepsilon > 0$,
    $ | \{ y\in G:  |h(g)- r|> \varepsilon \} | < \kappa$.
Equivalently,   for every $\varepsilon >0$,  there  exists   $B\in  [G]  ^{< \kappa} $   such that
$|h(x)- h(y)|< \varepsilon $  for all  $x, y\in  G\setminus B$.
Thus, if $h$ is not a convergent  at  infinity  then  there  exist  $\varepsilon > 0$ and the subsets
$ \{x_{\alpha}: \alpha < \kappa\}$, $\{y_{\alpha}: \alpha < \kappa\}$
  of  $G$  such that
  $|h(x_{\alpha})- h(y_{\alpha})|> \varepsilon $
   for each  $\alpha< \kappa$.

\vspace{7 mm}

{\bf Theorem 1}. {\it
If $\kappa< |G|$ then every  slowly  oscillating function
$h: (G, \mathcal{E}_{\kappa})\longrightarrow \mathbb{R}$
 is  convergent  at infinity, so the corona
 $\nu_{\kappa}(G)$ of $(G,\mathcal{E}_{\kappa})$
  is a  singleton.

\vspace{5 mm}

Proof.}
We need the following  auxiliary statement
\vspace{5 mm}

  $(\ast)$  {\it for every  subset  $A$ of  $G$  of  cardinality $\kappa$, there exists a  subgroup $S$  such that
   $A\subseteq  S$,  $|S|=\kappa $  and
   $h | _{SqS}$  is  constant  for each  $g\in  G \setminus S$.}
   \vspace{5 mm}

For  $\kappa = \omega$, $(\ast)$ follows  from Theorem 3.1(i)  in \cite{b4}.
We suppose  that   $\kappa > \omega$,  denote  by $cf \  \kappa$  the cofinality  of $\kappa$  and write   $A$  as  the
 union
 $\cup _{\alpha< cf \  \kappa} A_{\alpha}$,
  $|A_{\alpha}| < \alpha$
  for each
$   \alpha < cf  \  \kappa $.

We take the subgroup $S_{0}$  generated by $A_{0}$  and suppose that,  for some
 $\beta < cf \ \kappa $,  the  subgroups
 $\{S_{\alpha} : \alpha < \beta\}$  have been chosen.
  If  $\beta $  is a limit ordinal then we take the subgroup $S_{\beta}$ generated  by
  $A_{\beta}\cup \{S_{\alpha}: \alpha <\beta\}$.
If  $\beta$ is a non-limit ordinal then
$\beta = \gamma + n +1$,
  where $\gamma$ is a limit ordinal and $n\in \omega$.
Since $h$ is slowly oscillating, there exists
 $K\in [G] ^{<\kappa}  $  such that
 $$diam \  h(S_{\gamma+n}  \  g  \  S_{\gamma+n})< \frac{1}{n+1} $$
 for each $g\in G \setminus  K$.
We denote by  $S _{\gamma +n+1}$  the subgroup generated by
$ S _{\gamma+n} \cup A  _{\gamma+n+1} \bigcup K$.

After  $cf  \   \kappa $  steps we put
 $ S = \cup _{\alpha<cf \  \kappa}  S_{\alpha}  $.
By the construction,
$diam \  h(S g S ) < \frac{1}{n+1} $  for all  $n\in \omega$  and
$g\in G  \setminus  S$,   so  $(\ast)$  is proven.

We  suppose that $h$  is not convergent  at infinity and choose  $\varepsilon > 0$  and the subsets
 $\{ x_{\alpha} : \alpha <\kappa$,  $\{ y_{\alpha} : \alpha <\kappa\}$
  such  that
  $|h(x_{\alpha}) -  h(y_{\alpha} ) | > \varepsilon$   for  each $\alpha  <  \kappa$.
Then we put
$A=\{  x_{\alpha},  y_{\alpha} : \alpha < \kappa  \}$  and  choose the subgroup  $S$  satisfying  $(\ast)$.

We fix an arbitrary   $z\in G \setminus  S $.
Since  $h$  is slowly oscillating, there exists   $B\in [G] ^{<\kappa} $  such that
$$|h(zx)- h(x)|< \frac{\varepsilon}{2},  \   |h(zy)- h(y)|< \frac{\varepsilon}{2} $$
for all $x,y \in G\setminus B$.
Then we take
 $x_{\alpha}, y_{\alpha} \in G\setminus B$  and get
$$|h(z x_{\alpha}) - h(x_{\alpha})|< \frac{\varepsilon}{2}, \  |h(z y_{\alpha}) - h(y_{\alpha})|< \frac{\varepsilon}{2}. $$

By the choice of  $S, $ we have $x_{\alpha} \in  S$, $y_{\alpha} \in  S$ so
$h(z x_{\alpha}) = h(y_{\alpha} z)$
 and
$ |h( x_{\alpha}) - h(y_{\alpha})< \varepsilon$,
  a contradiction.    \hfill  $\Box$

 \vspace{8 mm}

For a cardinal $\kappa$,  we denote
$U_{\kappa}= \{ p\in \beta\kappa : |P| = \kappa$
 for each $P\in p\}$.

\vspace{7 mm}

{\bf Theorem 2}. {\it
If $\kappa=|G|$ and $\kappa$ is regular then
$\nu_{\kappa}(G)$
 contains a copy of  $U_{\kappa}$.

 \vspace{5 mm}

Proof.}
If $\kappa=\omega$ then, by [13, Theorem 2.1.1],  the ballean $(G, \mathcal{E}_{\omega})$ is metrizable and, by   [9, Theorem 1],  $\nu_{\omega}(G)$ contains a copy of $\beta\omega\setminus\omega$.

For $\kappa>\omega$,   we write  $G$ as the union of the increasing  chain
$\{G_{\alpha}: \alpha< \kappa\}$
 of subgroups such that
 $|G_{\alpha}|< \kappa$ for each $\alpha< \kappa$.
Since $\kappa$ is regular,  every function
$h:  (G, \mathcal{E}_{\kappa}) \longrightarrow \{0,1\}$
such that, for each
$\alpha< \kappa$,
$h|_{G_{\alpha+1}\setminus G_{\alpha}}= const $
 is slowly oscillating.
For every
$\alpha< \kappa$,
we pick
$x_{\alpha} \in G_{\alpha+1}\setminus G_{\alpha}$
 and put
 $X= \{x_{\alpha}: \alpha<\kappa  \}.$
We denote
$$U=\{p\in \beta G: X\in p, \   G\setminus G_{\alpha} \in p  \   \   for  \ \  each  \ \   \alpha < \kappa\}. $$
If
 $r, q\in U$ and $r\neq q$  then there is a slowly  oscillating function
  $h: (G, \mathcal{E}_{\kappa})\longrightarrow  \{0,1\} $  such that
  $h^{\beta}(r)\neq h^{\beta}(q).$
Clearly, $U$ is homeomorphic to $U_{\kappa} $.  \hfill  $\Box$

 \vspace{8 mm}

For a ballean  $(X, \mathcal{E})$, we  consider an
 equivalence
 $\sim_{\{0,1\}}$
 on $X^{\sharp}$
  defined by the rule:
  $p \sim_{\{0,1\}} q$
   if and only if
   $f^{\beta}(p) =  f^{\beta}(q)$
    for each  slowly oscillating function
    $f: X \longrightarrow\{0,1\}$.
The quotient
$X^{\sharp}/ \sim_{\{0,1\}}$
   is called
  \cite{b6}, \cite{b7}
     the {\it binary corona} of  $(X,  \mathcal{E})$.
For  a group  $G$,  the  binary  coronas of
$(G, \mathcal{R}_{\omega})$
 and
 $(G, \mathcal{E}_{\omega})$
 are known as the {\it  spaces of  ends and bi-ends}  of $G$ respectively \cite{b2}, \cite{b15}.

A  ballean  $(X, \mathcal{E})$ is called  {\it cellular} if  $\mathcal{E}$  has a base consisting of equivalence
relations.
If
$(X, \mathcal{E})$
 is cellular  then, by  [6, Lemma 4.3],
 $ \sim \  = \  \sim_{\{0,1\}}$
 so the  binary corona and corona of $(X, \mathcal{E})$  coincide.

If $\kappa > \omega$ then the coarse structure  $\mathcal{E}_{\kappa}$  of $G$  has  the base
$$\{\{ (x,y): y\in SxS\}: \  S\in [G]^{<\kappa} \ and \  S \
  is \ a \  subgroup \ of \  G\}. $$
It follows that the ballean
$(G, \mathcal{E}_{\kappa})$
 is cellular.

\vspace{7 mm}

{\bf Theorem 3 [11]}. {\it
Let $G$  be a group of  singular   cardinality $\kappa$. For any finite partition  $A_{1},\ldots  , A_{n}$  of  $G$,  there exist  $i\in \{ 1, \ldots  , n\}$  and  $S\in [G]  ^{<\kappa }$  such that } $G= S A_{i} S$.

\vspace{7 mm}

{\bf Theorem 4}. {\it
If  $\kappa= |G|$ and $\kappa$  is singular then the  corona $\nu_{\kappa}(G)$  is a singleton.

\vspace{5 mm}

Proof.}
Since $(G, \mathcal{E}_{\kappa})$ is cellular,  it  suffices to show that
$h^{\beta}(p)= h^{\beta}(q)$
 for any   $p,q\in G^{\sharp} $  and every
  slowly oscillating  function  $h: (G, \mathcal{E}_{\kappa})\longrightarrow\{0,1\} $.
We assume the contrary and choose  corresponding $p, q \in  G^{\sharp}$  and $h$.
Let
$h^{\beta} (p)= 1$,  $h^{\beta} (q)= 0$.
Then $G$  can be partitioned
$G=P \cup Q$  so  that
$h| _{P} = 1$, $h| _{Q} = 0$
 and $|P| = |Q|= \kappa$.
By Theorem 3, there exists
 $S\in [G]^{<\kappa} $  such that  either
 $G= SPS$  or $G= SQS$ ,  which is impossible because $h$ is slowly oscillating.  \hfill  $\Box$

\section{Comments}

1.	We suppose that  $\kappa< |G|$,  $A\subset G$, $|A|\geq \kappa$, $|G \setminus A|\geq\kappa$.
We define a function
$h: (G, \mathcal{E}_{\kappa})\longrightarrow\{0,1\}$
 by $h|_{A}= 1$, $h|_{G\setminus A}= 0$.
By Theorem 1,  $h$ is not slowly  oscillating.
Hence, there  exists $S\in [G] ^{<\kappa}  $  such that   $SAS\cap  (G \setminus A)\geq\kappa$.
Taras Banakh  noticed that, for $\kappa = \omega$,   this statement follows also from Theorem 1.20 in \cite{b2}.

\vspace{7 mm}

2.	With the same proof, Theorem 2  holds for the  corona  $\rho_{\kappa}(G)$
  of the ballean
  $(G, \mathcal{R}_{\kappa})$
   in place $\nu_{\kappa}(G)$.
Let $X$ be a set of cardinality  $|X|\geq \kappa$, $F(X) $  denotes  the free group in the alphabet $X$.
Example 3.3 from  \cite{b4}
and Theorem 4 from
\cite{b1}
 show that  $\rho_{\kappa}(F(X) )$ is not a  singleton.
Hence,  Theorem 1 and Theorem 4 do not  hold for  $G=F(X)$.
On the other hand, by  [4, Theorem 3.1],  the corona
$\rho_{\omega}(G)$
is a singleton provided that  $|G|>\omega $   and each  countable  subset of $G$ is contained  in some
countable normal subgroup.
\vspace{7 mm}

{\bf Question 1}.
{\it How can one detect  whether the corona  $\rho_{\omega}(G)$ is a singleton?}
\vspace{7 mm}

{\bf Question 2}.
{\it Are $\rho_{\omega}(G)$ and $\nu_{\omega}(G)$  homeomorphic for every countable group $G$?}

\vspace{7 mm}

By Theorem 2, for a countable group $G$,
$|\rho_{\omega}(G)|= \mid \nu _{\omega}(G)\mid  $
but the binary coronas of $(G, \mathcal{R}_{\omega} )$ and $(G, \mathcal{E}_{\omega} )$ could be of different cardinalities,  see  [2, Theorem 1.20].
\vspace{7 mm}

3. In \cite{b4},  \cite{b8} a function convergent at infinity is called constant at infinity.
This change in terminology  and reference \cite{b15}  are suggested by Yves  Cornulier.
 In  the case of an uncountable abelian group  $G$  of  cardinality $\kappa$,
   the infinity of the  space of
  $\kappa$-ends of $G$  is noticed in [15, Satz III].

\vspace{6 mm}

CONTACT INFORMATION

I.~Protasov: \\
Faculty of Computer Science and Cybernetics  \\
        Kyiv University  \\
         Academic Glushkov pr. 4d  \\
         03680 Kyiv, Ukraine \\ i.v.protasov@gmail.com

\end{document}